\documentclass[12pt]{amsart}
\pdfoutput=1
\usepackage{graphicx}
\usepackage{amssymb}
\usepackage{amsmath}
\usepackage{amsthm}
\usepackage{amscd}
\usepackage{epsfig}

\newtheorem{theorem}{Theorem}[section]

\newtheorem{lemma}[theorem]{Lemma}

\theoremstyle{definition}

\numberwithin{equation}{section}

\newcommand{\N}{\mathbb N}

\begin{document}

\title{A note on dimensional entropy for amenable group actions}

\author [Dou Dou and Ruifeng Zhang]{Dou Dou* and Ruifeng Zhang}

\address[D.~Dou]{Department of Mathematics, Nanjing University,
Nanjing, Jiangsu, 210093, P.R. China} \email{doumath@163.com}

\address[R.~Zhang]{School of Mathematics, Hefei University of Technology, Hefei, Anhui,
230009, P.R. China}
\email{rfzhang@mail.ustc.edu.cn}


\thanks{*Corresponding author}


\begin{abstract}
In this short note, for countably infinite  amenable group actions, we provide topological proofs for the following results: Bowen topological entropy (dimensional entropy) of the whole space equals the usual topological entropy along tempered F{\o}lner sequences; the Hausdorff dimension of an amenable subshift (for certain metric associated to some F{\o}lner sequence) equals its topological entropy. This answers questions by Zheng and Chen \cite{ZC} and Simpson \cite{S}.
\end{abstract}

\maketitle

%%%%%%%%%%%%%%%%%%%%%%%%%%%%%%%%%%%%%%%%%%%%%%%%%%%%%%%%%%%%%%%%%%%%%%%%%%%%%%%%%%%%%%%%%%%%%%%%%%%%%%%%%%%%%%%%%%%%%%%%%%%%%%%%%%%%
%%%%%%%%%%%%%%%%%%%%%%%%%%%%%%%%%%%%%%%%%%%%%%%%%%%%%%%        Introduction      %%%%%%%%%%%%%%%%%%%%%%%%%%%%%%%%%%%%%%%%%%%%%%%%%%%
%%%%%%%%%%%%%%%%%%%%%%%%%%%%%%%%%%%%%%%%%%%%%%%%%%%%%%%%%%%%%%%%%%%%%%%%%%%%%%%%%%%%%%%%%%%%%%%%%%%%%%%%%%%%%%%%%%%%%%%%%%%%%%%%%%%%
\section{Introduction}
Let $(X,G)$ be a $G-$action topological dynamical system, where $X$ is a compact Hausdorff space and $G$ a topological group. Throughout this paper, $G$ is always assumed to be a
countably infinite amenable group,
i.e. there exists a sequence of nonempty
finite subsets $\{F_n\}$ of $G$ (a {\it F{\o}lner sequence}) such that
$$\lim_{n\rightarrow+\infty}\frac{|F_n\vartriangle gF_n|}{|F_n|}=0, \text{ for all } g\in G.$$

Let $A$ and $K$ be two nonempty finite subsets of $G$. For $\delta>0$, the set $A$ is said to be {\it $(K, \delta)$-invariant} if
\begin{align*}
  \frac{|B(A,K)|}{|A|}<\delta,
\end{align*}
where $B(A,K)$, the {\it $K$-boundary} of $A$, is defined by
$$B(A,K)=\{g\in G: Kg\cap A\neq\emptyset \text { and } Kg\cap(G\setminus A)\neq\emptyset\}.$$
Another equivalent condition for the sequence of finite subsets $\{F_n\}$ of $G$ to be a F{\o}lner sequence is that $\{F_n\}$ becomes more and more
invariant, i.e. for any $\delta>0$ and any finite subset $K$ of $G$, $F_n$ is $(K,\delta)$-invariant
for sufficiently large $n$. For more information on amenable groups and their actions one may refer to \cite{C,KL,OW}.

For the case $G=\mathbb {Z}$, by resembling the definition of Hausdorff dimension, Bowen \cite{B} introduced a definition of topological entropy on subsets.
This definition is also known as dimensional entropy and has plenty of applications to thermodynamical formulism, fractal geometry, hyperbolic systems, multi-fractal analysis and so on (see, for example, \cite{FH} and \cite{P}).

For the case that $G$ is a general countably infinite amenable group, Bowen's dimensional entropy was recently introduced in \cite{ZC} in the following way.

Let $\mathcal{U}$ be a finite open cover of $X$. For a subset $F$ of $G$, denote by $\mathcal{W}_{F}(\mathcal{U})$ the
collection of families $\mathbf{U} = \{U_g\}_{g\in F}$ with $U_g\in \mathcal{U}$ (we also call $\mathbf{U}$ a $\mathcal{U}$-word or a $(\mathbf{U}, F)$-name). For $\mathbf{U}\in\mathcal{W}_{F}(\mathcal{U})$ we call the integer
$m(\mathbf{U}) = |F|$ the length of $\mathbf{U}$ and the set $F$ the domain of $\mathbf{U}$ (denoted by ${\rm dom}(\mathbf{U})$). For any $\mathcal{U}$-word $\mathbf{U}$, define
\begin{align*}
\mathbf{X}(\mathbf{U}) &:= \bigcap_{g\in {\rm dom}(\mathbf{U})} g^{-1}U_g \\
&=\{x\in X : gx \in U_g \text{ for }g\in {\rm dom}(\mathbf{U})\}.
\end{align*}

Now let $\{F_n\}$ be a F{\o}lner sequence in $G$.
For $Z\subset X$, we say that a collection of $(\mathbf{U}, F_n)$-names $\Lambda\subset \bigcup_{n\ge 1}\mathcal{W}_{F_n}(\mathcal{U})$ covers $Z$ if
$\bigcup_{\mathbf{U}\in\Lambda}\mathbf{X}(\mathbf{U})\supset Z$. For $s\ge 0$,
define
$$\mathcal{M}(Z,\mathcal{U}, N,s,\{F_n\})=\inf\limits_{\Lambda}\{\sum\limits_{\mathbf{U}\in \Lambda}\exp(-s m(\mathbf{U}))\},
$$
where the infimum is taken over all $\Lambda\subset \bigcup_{j\ge N}\mathcal{W}_{F_j}(\mathcal{U})$ that covers $Z$.
It is not hard to see that $\mathcal{M}(\cdot,\mathcal{U}, N,s,\{F_n\})$ is a finite outer measure on $X$.
As $\mathcal{M}(Z,\mathcal{U}, N,s,\{F_n\})$ increases when $N$ increases, then define
$$\mathcal{M}(Z,\mathcal{U},s,\{F_n\})= \lim\limits_{N\rightarrow +\infty}\mathcal{M}(Z,\mathcal{U}, N,s,\{F_n\})$$
and
\begin{align*}
    h^B_{top}(\mathcal{U}, Z,\{F_n\}) &= \inf\{s\ge 0 : \mathcal{M}(Z,\mathcal{U},s,\{F_n\})= 0\} \\
    &= \sup\{s\ge 0 : \mathcal{M}(Z,\mathcal{U},s,\{F_n\}) = +\infty\}.
\end{align*}
Set
$$h^B_{top}(Z,\{F_n\})= \sup_{\mathcal {U}}
h^B_{top}(\mathcal{U}, Z,\{F_n\}),$$
where $\mathcal{U}$ runs over finite open covers of $X$. We call $h^B_{top}(Z,\{F_n\})$ the Bowen topological
entropy or dimensional entropy of $Z$(w.r.t. the F{\o}lner sequence $\{F_n\}$).

Recall that a F{\o}lner sequence $\{F_n\}$ in $G$ is said to be {\it tempered} if there exists a constant $C>0$ which is independent of $n$ such that
\begin{align}\label{tempered}
|\bigcup_{k<n}F_k^{-1}F_n|\le C|F_n|, \text{ for any }n.
\end{align}

In \cite{ZC}, Zheng and Chen proved that for any tempered F{\o}lner sequence $\{F_n\}$ in $G$
with the increasing condition $\lim\limits_{n\rightarrow+\infty}\frac{|F_n|}{\log n}=\infty$,
the Bowen topological entropy $h_{top}^B(X,\{F_n\})$ coincides with the usual topological entropy.
This generalized Bowen's classical result in \cite{B} to discrete countable amenable group actions.

Differs from Bowen's original topological proof, Zheng and Chen used tools in ergodic theory---they employed SMB theorem, Brin-Katok's formula for local entropy and variation principle for topological entropy. They also asked in \cite{ZC} whether there exists a pure topological proof.

In this note, we will give a direct proof for this result. More precisely, we will prove the following
\begin{theorem}[Main theorem]\label{th-1-1}
Let $(X,G)$ be a $G-$action topological dynamical system where $X$ is a compact hausdorff space and $G$ is a discrete countable amenable group, then for any tempered F{\o}lner sequence $\{F_n\}$ in $G$,
we have
$$h_{top}^B(X,\{F_n\})=h_{top}(X,G),$$
where $h_{top}(X,G)$ is the topological entropy of $(X,G)$ defined through open covers.
\end{theorem}

The itinerary of our proof uses Simpson's idea for proving the relations between entropy and Hausdorff dimension of $\mathbf{Z}^d$-subshifts (the idea was originally
by Furstenberg \cite{F} and it was also used in Bowen \cite{B}). But there exists essential difficulty for amenable group actions.
The key point is that we need a Vitali type covering lemma. To this aim we employ
a more general covering lemma developed by Lindenstrauss \cite{L} (for proving the amenable version of the pointwise ergodic theorem).
In section 2 we will give the detailed proof.

\section{Proof of the main theorem}
We first recall the covering lemma by Lindenstrauss.

\begin{lemma}[Corollary 2.7 of Lindenstrauss \cite{L}]\label{lemma-2-1}
  For any $\delta\in (0,1/100), C>0$ and finite $D\subset G$, let $M\in\N$ be sufficiently large (depending only on $\delta,C$ and $D$).
  Let $F_{i,j}$ be an array of finite subsets of $G$ where $i=1,\ldots,M$ and $j=1,\ldots,N_i$, such that
  \item [1.] For every $i$, $\bar{F}_{i,*}=\{F_{i,j}\}_{j=1}^{N_i}$ satisfies
  $$|\bigcup_{k'<k}F_{i,k'}^{-1}F_{i,k}|\le C|F_{i,k}|,\quad \text{ for }k=2,\ldots,N_i.$$
  \quad Denote $F_{i,*}=\cup\bar{F}_{i,*}$.
  \item [2.] The finite set sequences $F_{i,*}$ satisfy that for every $1<i\le M$ and every $1\le k \le N_i$,
  $$|\bigcup_{i'<i}DF_{i',*}^{-1}F_{i,k}|\le (1+\delta)|F_{i,k}|.$$

  Assume that $A_{i,j}$ is another array of finite subsets of $G$ with $F_{i,j}A_{i,j}\subset F$ for some finite subset $F$ of $G$.
  Let $A_{i,*}=\cup_jA_{i,j}$ and $$\alpha=\frac{\min_{1\le i\le M}|DA_{i,*}|}{|F|}.$$

Then the collection of subsets of $F$,
  $$\tilde{\mathcal{F}}=\{F_{i,j}a: 1\le i\le M, 1\le j\le N_i \text{ and }a\in A_{i,j}\}$$
  has a subcollection $\mathcal{F}$ that is $10\delta^{1/4}$-disjoint such that
  $$|\cup \mathcal{F}|\ge (\alpha-\delta^{1/4})|F|.$$
\end{lemma}

We note here that a collection $\tilde{\mathcal{F}}$ of finite subsets of $G$ is said to be $\delta$-disjoint if for every $A\in \tilde{\mathcal{F}}$
there exists an $A'\subset F$ such that $|A'|\ge (1-\delta)|A|$ and such that $A\cap B=\emptyset$ for every $A\neq B\in \tilde{\mathcal{F}}$.

Now we turn to prove our main theorem. The proof of the upper bound is straight-forward and has been shown in Section 4 of \cite{ZC}. We omit it here.

In the following we will give the proof of the lower bound.

For any $\epsilon>0$, let $0<\delta<\min\{\epsilon,1/100\}$ be small enough such that
\begin{align}\label{condition-2-1}
  -(1-2\delta-11\delta^{1/4})\log(1-2\delta-11\delta^{1/4})-(2\delta+11\delta^{1/4})\log(2\delta+11\delta^{1/4})<\epsilon.
\end{align}
Let $D=\{e_G\}\subset G$ and let $C>0$ be the constant in the tempered condition \eqref{tempered} for the F{\o}lner sequence $\{F_n\}$.
Let $M>0$ be large enough to satisfy the requirement of Lemma \ref{lemma-2-1} corresponding to $\delta, D$ and $C$.

Let $\mathcal{U}$ be any finite open cover of $X$ and $s>0$ such that
$\mathcal{M}(X,\mathcal{U},s,\{F_n\})=0$. We will show $h_{top}(G,\mathcal{U})\le s$ and then it follows that
$h_{top}^B(X,\mathcal{U},\{F_n\})\ge h_{top}(G,\mathcal{U})$.

Since $\mathcal{M}(X,\mathcal{U},s,\{F_n\})=\lim\limits_{N\rightarrow +\infty}\mathcal{M}(X,\mathcal{U}, N,s,\{F_n\})=0$, for each $i=1,2,\ldots$,
there exist $0<p_i$ and $\Lambda_i\subset \bigcup_{j\ge p_i}\mathcal{W}_{F_j}(\mathcal{U})$ such that $\Lambda_i$ covers $X$ and
$$\sum\limits_{\mathbf{U}\in \Lambda_i}\exp(-s m(\mathbf{U}))< 2^{-i}.$$ As $X$ is compact, we may let each $\Lambda_i$ be finite.
In addition, we let the sequence $\{p_i\}$ be increasing and
\begin{align}\label{condition-2-2}
  |F_j|>\frac{1}{\delta(1-10\delta^{1/4})}, \text{ for all }j\ge p_1.
\end{align}

Denote $\Lambda_{\infty}=\bigcup\limits_{i=1}^{\infty}\Lambda_i$. Then
$$\sum_{\mathbf{U}\in \Lambda_{\infty}}\exp(-s m(\mathbf{U}))< \sum_{i=1}^{\infty}2^{-i}=1.$$
Hence
\begin{align*}
  &\sum_{k=1}^{\infty}\sum_{\mathbf{U}_1,\ldots,\mathbf{U}_k\in \Lambda_{\infty}}\exp\big(-s (m(\mathbf{U}_1)+\ldots+m(\mathbf{U}_k))\big)\\
  =&\sum_{k=1}^{\infty}\bigg (\sum_{\mathbf{U}\in \Lambda_{\infty}}\exp\big(-s m(\mathbf{U})\big)\bigg )^k\\
  =&S<\infty.
\end{align*}

For each $i=1,\ldots, M$, let $\{F_{n_{i,1}},\ldots,F_{n_{i,N_i}}\}=\{{\rm dom}(\mathbf{U}): \mathbf{U}\in \Lambda_i\}$ with $n_{i,1}<\ldots<n_{i,N_i}$. Then let $\{F_{i,1},\ldots,F_{i,N_i}\}$ in Lemma \ref{lemma-2-1} be as
$$\{F_{i,1},\ldots,F_{i,N_i}\}:=\{F_{n_{i,1}},\ldots,F_{n_{i,N_i}}\}.$$

For any $x\in X$ and sufficiently large $n$ (independent on $x$), let
\begin{align*}
  A_{i,j}=&\{a\in F_n: F_{i,j}a\subset F_n \text{ and there exists }\mathbf{U}\in \Lambda_i \\
  &\qquad \qquad \text{ such that } {\rm dom}(\mathbf{U})=F_{i,j} \text{ and }ax\in \mathbf{X}(\mathbf{U})\}.
\end{align*}
We note here that $A_{i,j}$ depends on $x$. Denote by $e_G$ the identity element of $G$. For any $g\in F_n\setminus B(F_n,F_{i,*}\cup\{e_G\})$, we have $F_{i,*}g\subset F_n$.
Hence for all $1\le j\le N_i$, $F_{i,j}g\subset F_n$. Since $\Lambda_i$ covers $X$, there exists some $\mathbf{U}\in \Lambda_i$ which contains $gx$.
Suppose ${\rm dom}(\mathbf{U})=F_{i,j}$. We then have $g\in A_{i,j}$. This implies that
$$A_{i,*}\supset F_n\setminus B(F_n,F_{i,*}\cup\{e_G\}).$$
Hence when $n$ is large enough such that $F_n$ is
$(F_{i,*}\cup\{e_G\},\delta)$-invariant for all $1\le i\le M$ and $1\le j\le N_i$,
$$\alpha=\frac{\min_{1\le i\le M}|DA_{i,*}|}{|F_n|}>1-\delta.$$

The first requirement for the array $\{F_{i,j}\}$ in Lemma \ref{lemma-2-1} comes from the tempered condition of $\{F_n\}$. To ensure the second requirement,
we can make $p_i$ be large enough compared with $n_{i-1,N_{i-1}}$ for every $2\le i\le M$.

Now we can apply Lemma \ref{lemma-2-1} to the collection of subsets of $F_n$,
  $$\tilde{\mathcal{F}}=\{F_{i,j}a: 1\le i\le M, 1\le j\le N_i \text{ and }a\in A_{i,j}\}$$
to obtain a subcollection $\mathcal{F}$ that is $10\delta^{1/4}$-disjoint such that
  $$|\cup \mathcal{F}|\ge (1-\delta-\delta^{1/4})|F_n|.$$

For each element in $\mathcal F$, we will fix one way to write it as the form $F_{i,j}a$.
Denote by $\bar{A}$ the collection of $a$'s such that $F_{i,j}a$ occurs in $\mathcal{F}$. The cardinality of $\bar{A}$ is no more than the cardinality of the subcollection
$\mathcal{F}$. Since $\mathcal{F}$ is $10\delta^{1/4}$-disjoint, we have that
$$\sum_{F_{i,j}a\in \mathcal{F}}|F_{i,j}a|\le \frac{1}{1-10\delta^{1/4}}|\cup \mathcal{F}|\le \frac{1}{1-10\delta^{1/4}}|F_n|.$$
Hence
\begin{align*}
  |\bar{A}|\le |\mathcal{F}|\le \frac{1}{\min |F_{i,j}|}\cdot\frac{1}{1-10\delta^{1/4}}|F_n|\le \delta|F_n| \quad (\text{by \eqref{condition-2-2}}).
\end{align*}
We note that for each $a\in \bar{A}$, there may exist more than one $F_{i,j}$'s with $F_{i,j}a\in \mathcal{F}$. Denote by $n(a)$ the total number of such $F_{i,j}$'s.

We now construct $\mathbf{U}(x)=\{U_g(x)\}_{g\in F_n}\in \mathcal{W}_{F_n}(\mathcal{U})$, a $(\mathbf{U}, F_n)$-name associated to $x$, in the following way:
\begin{enumerate}
  \item [1.] For each $g$ that lies in exactly one element of $\mathcal{F}$, say $F_{i,j}a$, assume that $ax\in \mathbf{X}(\mathbf{U})$ for $\mathbf{U}=\{U_g\}_{g\in F_{i,j}}\in \Lambda_i$.
We then let $U_g(x)=U_{ga^{-1}}$.
\item [2.] For each $g$ that lies in either $F_n\setminus \bigcup\mathcal{F}$ or the overlapping part of elements in $\mathcal{F}$, we choose any $U\in\mathcal{U}$ that contains $gx$
and then let $U_g(x)=U$.
\end{enumerate}
Clearly $x\in \mathbf{X}(\mathbf{U}(x))\in \mathcal{U}^{F_n}$ and hence $\{\mathbf{X}(\mathbf{U}(x))\}_{x\in X}$ forms a subcover of $\mathcal{U}^{F_n}$.

In the following we will estimate the cardinality of the subcover $\{\mathbf{X}(\mathbf{U}(x))\}_{x\in X}$.
Let $\Lambda=\bigcup_{x\in X}\{\mathbf{U}(x)\}$. An element $\mathbf{U}=\{U_g\}_{g\in F_n}$ in $\Lambda$ is determined by two parts:
\begin{enumerate}
\item [1.] the $U_g$'s for the $g$'s that lie in exactly one element of $\mathcal{F}$;
\item [2.] the $U_g$'s for the $g$'s that lies in either $F_n\setminus \bigcup\mathcal{F}$ or the overlapping part of elements in $\mathcal{F}$.
\end{enumerate}
The first part corresponds to the choice of $\mathcal{F}$ and the choice of the sequence of $\mathbf{U}_1, \ldots, \mathbf{U}_k$ with $m(\mathbf{U}_1)+\ldots+m(\mathbf{U}_k)\le \frac{1}{1-10\delta^{1/4}}|F_n|$ (here $k$ is the cardinality of $\mathcal{F}$). When $\mathbf{U}_1, \ldots, \mathbf{U}_k$ are fixed, if we know $\bar{A}$ and
the number $n(a)$ for each $a\in \bar{A}$, $\mathcal{F}$ is then determined.
Hence the first part of $\mathbf{U}\in \Lambda$ can be determined by $\bar{A}$, $\{n(a)\}_{a\in \bar{A}}$ with $\sum_{a\in \bar{A}}n(a)=k$ and the sequence of $\mathbf{U}_1, \ldots, \mathbf{U}_k$ with $m(\mathbf{U}_1)+\ldots+m(\mathbf{U}_k)\le \frac{1}{1-10\delta^{1/4}}|F_n|$.
The second part corresponds to the choice of $U_g$'s for $g$ in either $F_n\setminus \bigcup\mathcal{F}$ or the overlapping part of $\bigcup\mathcal{F}$.

Since $|\bar{A}|\le \delta|F_n|$, the total number of the choices of $\bar{A}$'s is bounded from above by
$\sum_{m=1}^{\lfloor \delta|F_n|\rfloor}\binom{|F_n|}{m}$. Applying the Stirling formula
\begin{align*}
  n!=\sqrt{2\pi n}(\frac{n}{e})^ne^{\alpha_n}, \frac{1}{12n+1}<\alpha_n<\frac{1}{12n},
\end{align*}
we have that
\begin{align*}
 \sum_{m=1}^{\lfloor \delta|F_n|\rfloor}\binom{|F_n|}{m}&\le\delta|F_n|\cdot\frac{\sqrt{2\pi |F_n|}e^{\alpha_{|F_n|}-\alpha_{|F_n|-\lfloor\delta|F_n|\rfloor}-\alpha_{\lfloor\delta|F_n|\rfloor}}}{\sqrt{2\pi (|F_n|-\lfloor\delta|F_n|\rfloor)}\sqrt{2\pi \lfloor\delta|F_n|\rfloor}}\cdot\\ &\ \ \ \ \ \ \ \ \quad\quad (\frac{|F_n|}{|F_n|-\lfloor\delta|F_n|\rfloor})^{|F_n|-\lfloor\delta|F_n|\rfloor}\cdot(\frac{|F_n|}{\lfloor\delta|F_n|\rfloor})^{\lfloor\delta|F_n|\rfloor}\\
 &\le\delta|F_n|\cdot\frac{\sqrt{2\pi |F_n|}e^{\alpha_{|F_n|}-\alpha_{|F_n|-\lfloor\delta|F_n|\rfloor}-\alpha_{\lfloor\delta|F_n|\rfloor}}}{\sqrt{2\pi (|F_n|-\lfloor\delta|F_n|\rfloor)}\sqrt{2\pi \lfloor\delta|F_n|\rfloor}}\cdot\\ &\ \ \ \ \ \ \ \ \quad\quad (\frac{1}{1-\delta})^{(1-\delta)|F_n|+1}(\frac{1}{\delta})^{\delta|F_n|}\cdot(\frac{\delta|F_n|}{\lfloor\delta|F_n|\rfloor})^{\lfloor\delta|F_n|\rfloor}\\
 &\le\delta|F_n|\cdot\frac{\sqrt{2\pi |F_n|}e^{\alpha_{|F_n|}-\alpha_{|F_n|-\lfloor\delta|F_n|\rfloor}-\alpha_{\lfloor\delta|F_n|\rfloor}}}{\sqrt{2\pi (|F_n|-\lfloor\delta|F_n|\rfloor)}\sqrt{2\pi \lfloor\delta|F_n|\rfloor}}\cdot \\ &\ \ \ \ \ \ \ \ \quad \quad
 \frac{1}{1-\delta}\cdot e \cdot (\frac{1}{1-\delta})^{(1-\delta)|F_n|}\cdot(\frac{1}{\delta})^{\delta|F_n|}\\
 &=Q(|F_n|)\cdot \exp\bigg(\big(-(1-\delta)\log(1-\delta)-\delta\log \delta\big)|F_n|\bigg),
\end{align*}
where
\begin{align*}
  Q(n)=\frac{\delta}{1-\delta} n\cdot \sqrt{\frac{n}{2\pi (n-\lfloor\delta n\rfloor)\lfloor\delta n\rfloor}}\cdot e^{\alpha_n-\alpha_{n-\lfloor\delta n\rfloor}-\alpha_{\lfloor\delta n\rfloor}+1}
\end{align*}
We can choose $n$ sufficiently large to make $Q(|F_n|)\le \exp(\epsilon|F_n|)$. Together with \eqref{condition-2-1}, we can obtain that the total number of the choices of $\bar{A}$'s will not exceed $\exp(2\epsilon|F_n|)$ when $n$ is large enough.

The total number of the choices of the sequence $\{n(a)\}_{a\in\bar{A}}$ is no larger than $\binom{|\mathcal{F}|}{|\bar{A}|}$, which can be bounded by $\exp(\epsilon|F_n|)$.

 The cardinality of the union of $F_n\setminus \bigcup\mathcal{F}$ and the overlapping part of $\bigcup\mathcal{F}$ is no more than $(\delta+\delta^{1/4}+10\delta^{1/4})|F_n|$.
 Also note that for each $g$ of such case, $U_g$ has $\#\mathcal{U}$ many choices. Hence the cardinality of $\Lambda$
is no more than $\exp((\delta+11\delta^{1/4})|F_n|\log\#\mathcal{U})\cdot \exp(2\epsilon|F_n|)\cdot\exp(\epsilon|F_n|)$ times the total number of the sequences of $\mathbf{U}_1, \ldots, \mathbf{U}_k\in \Lambda_{\infty}$ with $m(\mathbf{U}_1)+\ldots+m(\mathbf{U}_k)\le \frac{1}{1-10\delta^{1/4}}|F_n|$.

Finally,
\begin{align*}
  N(\mathcal{U}^{F_n})2^{-s|F_n|/(1-10\delta^{1/4})}\le & \exp((\delta+11\delta^{1/4})|F_n|\log\#\mathcal{U})\cdot \exp(3\epsilon|F_n|)\\
  &\qquad \cdot\sum_{m(\mathbf{U}_1)+\ldots+m(\mathbf{U}_k)\le \frac{1}{1-10\delta^{1/4}}|F_n|}2^{-s|F_n|/(1-10\delta^{1/4})}\\
  \le & \exp\big(((\delta+11\delta^{1/4})\log\#\mathcal{U}+3\epsilon)|F_n|\big)\\
  &\qquad \cdot\sum_{\mathbf{U}_1, \ldots, \mathbf{U}_k\in \Lambda_{\infty}}2^{-s(m(\mathbf{U}_1)+\ldots+m(\mathbf{U}_k))}\\
  \le & \exp\big(((\delta+11\delta^{1/4})\log\#\mathcal{U}+3\epsilon)|F_n|\big)S.
\end{align*}
Letting $\epsilon$ tend to $0$, it holds that $h_{top}(G,\mathcal{U})\le s$.

\subsection*{Remark}
In \cite{S}, Simpson discussed the relation between the entropy and Hausdorff dimension for $\mathbf{Z}^d$ subshifts and showed that they are the same for some ``standard
metric". He also asked (Question 5.4 of \cite{S}, in part) whether it also holds for wider classes of groups or semigroups such as an amenable group. We remark here that the same argument
can be applied to show that the entropy and Hausdorff dimension for amenable subshifts (for some metric associated with some tempered F{\o}lner sequence) coincide. We can also show this by applying Theorem \ref{th-1-1}.

Recall that when $X$ is a metric space with metric $d$, the dimensional entropy for the system $(X,G)$ can be defined in
the following alternative way \cite{ZC}.

For a finite subset $F$ in $G$, denote $B_{F}(x,\epsilon):=\{y\in X: d(gx,gy)<\epsilon, \text{ for any }g\in F\}$.
For $Z\subseteq X, s\ge0, N\in\mathbf{N}$, $\{F_n\}$ a F{\o}lner sequence in $G$ and $\epsilon> 0$, define
$$\mathcal{M}(Z,N,\epsilon,s,\{F_n\}) = \inf\sum_i\exp(-s|F_{n_i}|),$$
where the infimum is taken over all countable families $\{B_{F_{n_i}}(x_i,\epsilon)\}$ such that
$x_i\in X, n_i\ge N$ and $\bigcup_i B_{F_{n_i}}(x_i,\epsilon)\supseteq Z$. Then let
$$\mathcal{M}(Z,\epsilon,s,\{F_n\})=\lim_{N\rightarrow+\infty}\mathcal{M}(Z,N,\epsilon,s,\{F_n\}),\mathcal{M}(Z,s,\{F_n\})=\lim_{\epsilon\rightarrow0}\mathcal{M}(Z,\epsilon,s,\{F_n\}).$$
Bowen topological entropy $h^B_{top}(Z,\{F_n\})$ can be equivalently defined as the critical value
of the parameter $s$, where $\mathcal{M}(Z,s,\{F_n\})$ jumps from $+\infty$ to $0$.

Let $A$ be a finite set. Consider the left action of $G$ on the product space $A^G=\{(x_g)_{g\in G}: x_g\in A\}$:
$$g'(x_g)_{g\in G}=(x_{gg'})_{g\in G}, \text{ for all }g'\in G\text{ and }(x_g)_{g\in G}\in A^G.$$
Let $X$ be a closed $G$-invariant subset of $A^G$ (we call such $X$ a {\it subshift}). Fix any tempered F{\o}lner sequence $\{F_n\}$ of $G$ with $F_1\subseteq F_2\subseteq\ldots$
and $\bigcup_{n}F_n=G$.
We can then define a metric $d$ on $A^G$ by the following: $d(x,y)=1$ if $x$ and $y$ are not equal on $F_1$; otherwise $d(x,y)=e^{-|F_n|}$ where $n$ is the maximal number such that $x_g=y_g$ for all $g\in F_n$. The Hausdorff measure
of $X$ is defined as
$$\mu_s(X)=\lim_{\epsilon\rightarrow 0}\inf_{\mathcal{E}}\sum_{E\in\mathcal{E}}{\rm diam}(E)^s,$$
where $\mathcal{E}$ takes over (open) coverings of $X$ such that ${\rm diam}(E)\le \epsilon$ for all $E\in\mathcal{E}$.
It is easy to see that the covering $\mathcal{E}$ can be also taken from the cylinders of $A^G$ with the form $[\omega]_{F_n}:=\{x\in A^G: x_g=\omega_g, \text{ for every }g\in F_n\}$, $\omega\in A^G$.
Hence
$$\mu_s(X)=\lim_{N\rightarrow \infty}\inf_{\mathcal{E}}\sum_{E\in\mathcal{E}}{\rm diam}(E)^s,$$
where $\mathcal{E}=\{[\omega^{i}]_{F_{n_i}}\}$ is taken over countable families of cylinders that cover $X$ with $n_i\ge N$.
The Hausdorff dimension of $X$ is then defined by
$$\dim(X)=\inf\{s\ge 0:\mu_s(X)=0\}=\sup\{s\ge 0:\mu_s(X)=\infty\}.$$

It is easy to see that $\dim(X)=h^B_{top}(X,\{F_n\})$ by comparing the definition of the Bowen topological entropy and Hausdorff dimension. By Theorem \ref{th-1-1},
we have that $\dim(X)=h_{top}(X,G)$ for the metric $d$ associated with any tempered F{\o}lner sequence $\{F_n\}$ of $G$ with $F_1\subseteq F_2\subseteq\ldots$
and $\bigcup_{n}F_n=G$.

{\bf Acknowledgements}
The authors would like to express their gratitude to the referees for their valuable suggestions and comments.
This work was done when both of the authors were visiting scholars at SUNY at Buffalo and they would like to thank CSC and SUNY at Buffalo for the support.
The first author was also partially supported by NNSF of China (Grant No. 11431012, 11401220). The second author was also
partially supported by NNSF of China (Grant No. 11671094).

\end{document}